\documentclass[12pt]{article}

\usepackage[a4paper]{geometry}
\usepackage{ajc-FR}
\usepackage{amsmath,amsfonts,latexsym}
\usepackage{graphicx}
\usepackage{hyperref}
\usepackage{placeins}

\Volume{XY(Z)}
\Year{2021}
\firstpage{1}
\lastpage{11}
\Revised{24 June 2021}
\runninghead{\copyright \, F. Rowley 2021\,/\, GLRGC FINAL DRAFT arxiv version, pages 1 -- 11}

\numberwithin{equation}{section}

\newtheorem{theorem}{Theorem}[section]

\newtheorem{definition}[theorem]{Definition}

\newenvironment*{proof}
{\begin{list}{}{\setlength{\leftmargin}{0em}\setlength{\rightmargin}{0em}}
\item[] {\sc Proof:}} {\hfill$\Box$
\end{list}}

\begin{document}

\title{\bf A generalised linear Ramsey graph construction}
\author{Fred Rowley\thanks{formerly of Lincoln College, Oxford, UK.}}
\Addr{West Pennant Hills, \\NSW, Australia.\\{\tt fred.rowley@ozemail.com.au}}

\date{\dateline{6 June 2019}{DD Mmmmm CCYY}}

\maketitle


\begin{abstract}

A construction described by the current author in 2017 uses two general linear graphs as `prototypes' to build a linear compound graph with inherited Ramsey properties.  

This paper describes a generalisation of that construction which has produced improved lower bounds for multicolour Ramsey numbers in many cases.  The resulting compound graphs are linear, as before, and under certain particular conditions, they can be cyclic.  

The mechanism of the new construction requires that the first prototype contains a triangle-free `template', with defined properties, in one colour.  This paper shows that in the compound graph, clique numbers in the colours of the first prototype may exceed those of the prototype.  However, it proves necessary only to test for a limited subset of the possible cliques using these colours, in order to evaluate the relevant clique numbers for the entire graph.  Clique numbers in the colours of the second prototype are equal to those of that prototype.  It has also been found that there are a number of useful cases in which the clique numbers of the first prototype are not increased by the compounding process.  These attributes enable the efficient searching of a new family of graphs of unlimited size. 

As a result of this construction many lower bounds can be improved.  The improvements include $R_4(5) \ge 4073$, $R_5(5) \ge 38914$, $R_3(6) \ge 1106$, $R_4(6) \ge 21302$, $R_4(7) \ge 84623$ and $R_3(9) \ge 14034$.  Furthermore, it is shown that $R_3(9) \ge 14081$ using a non-linear construction.   

The construction also enables improved lower bounds on $\lim_{\substack{r \rightarrow \infty}} R{_r}(k)^{1/r}$ for many values of $k$.

\end{abstract}

\bigskip

\bigskip

\section{Introduction}
This paper addresses the properties of undirected loopless graphs with edge colourings in an arbitrary number of colours, and the corresponding multicolour classical Ramsey numbers, by describing a powerful construction for multicolour Ramsey graphs. 

We first define a {\it Ramsey graph} $U(k_1, {\dots} \, ,k_r; m)$, with all $k_s \ge 2$, as a complete graph of order $m$ with a colouring such that for each colour $s$, where $1 \leq s \leq r$, there exists no monochrome copy of a complete graph $K_t$ which is a subgraph of $U$ in the colour $s$ for any $t \geq k_s$.  Equivalently, the maximum order of any such copy of $K_t$ in any colour $s$ is strictly less than $k_s$.  Such a graph $U$ may conveniently be described as a $(k_1, {\dots} , k_r; m)${\it-graph}.  

The {\it Ramsey number} $R(k_1, \dots \, , k_r)$ is the unique lowest integer $m$ such that no \newline
 $U(k_1, {\dots} \, , k_r; m)$ exists, and its existence is proved by Ramsey's Theorem.  When all the $k_i$ are equal, this may be written $R_r(k)$ and is referred to as a {\it diagonal} Ramsey number.  

The construction described in \cite{Rowley1} extended a previous construction of Giraud \cite{Giraud1}.  It allows the creation of linear Ramsey graphs with specific properties by combining the distance sets of two linear graphs, which we will refer to as {\it prototypes}.  We may refer to graphs constructed in this way generically as {\it compound graphs}.

In \cite{Giraud1}, by identifying symmetrical Schur partitions with the corresponding 'regular colourings' of Kalbfleisch \cite{Kalb}, Giraud effectively established that various infinite series of Ramsey graphs could be constructed, adding one colour at a time, and that the orders of those graphs could be determined by a simple formula based on the number of colours.  His results set lower bounds for a wide range of multicolour Ramsey numbers $R_r(k)$ of unlimited size, and for any $k$.  

In the diagonal case, the growth rate of the orders of the graphs in such series can often be seen to approach a fixed limit as the number of colours tends to infinity, which we will refer to as the {\it limiting growth rate} of the series.  Where the series arises from the repeated use of the same construction with the same prototype graph, we may associate the limiting growth rate with that combination.  

Giraud's results effectively established lower bounds on $\lim_{\substack{r \rightarrow \infty}} R{_r}(k)^{1/r}$ which for convenience we may write as $\Gamma(k)$.  The existence of $\Gamma(k)$ had been proved already in \cite{AbbM}.

For $k = 3$, we note here that Giraud's work was strongly foreshadowed, many years before, by a closely related result of Schur published in \cite{Sch1}.  In fact, for the specific case $k = 3$, it had already been somewhat bettered by Abbott and Moser in \cite{AbbM} who proved (after adjusting terminology) that $R{_r}(3) \ge 89^{(r/4-c.log(r))}$ for some fixed constant $c$, from which it follows that $\Gamma(3) \ge 3.071\dots$.  Much later, and after other significant work including \cite{AbbH}, it was proved in \cite{XXER} that $\Gamma(3) \ge 3.199\dots$.   

Building on this background, section 3 of this paper describes a generalisation of the construction demonstrated in \cite{Rowley1}, which produces superior lower bounds in many cases.  It can be applied to produce diagonal or non-diagonal Ramsey graphs.  The resulting graphs are always linear, and under particular conditions, they are cyclic.

This in turn allows us to extend the construction of infinite series of Ramsey graphs, aiming in particular to throw more light on lower bounds for $\Gamma(k)$ for any value of $k$.  







We also establish improved lower bounds on a range of `small' diagonal multicolour Ramsey numbers.  Section 4 records some selected computational results for individual Ramsey graphs.  These establish improved lower bounds for a range of graphs, including $R_4(5)$ and $R_5(5)$, $R_3(6)$ and $R_4(6)$, $R_4(7)$, and $R_3(9)$.  

In section 5, the key Tables from \cite{Rowley1} are updated to reflect the limiting growth rates implied by these improved results in particular cases for the previous construction.  Lower bounds for some specific $\Gamma(k)$ are derived by taking selected maxima of these growth rates, or extended using the new construction.  

Section 6 draws some conclusions from this work.

We begin with some further definitions of required notation in section 2.  

\section{Definitions and Notation}

In this paper, 

$K_n$ denotes the complete graph with order $n$.  

If $U$ denotes a complete graph $K_m$ with $m$ vertices $\{u_0, {\dots} , u_{m-1}\}$, then: 

A {\it (q-)colouring} of $U$ is a mapping of the edges $({u_i}, {u_j})$ of $U$ into the set of integers $s$ where $1\le s \le q$.  

The {\it length} of the edge $({u_i}, {u_j})$ is defined as $\mid j - i \mid$.  A length is often referred to as an {\it edge-length} in this paper, for clarity.

A colouring of $U$ is {\it linear} if and only if the colour of any edge $({u_i}, {u_j})$ depends only on the length of that edge.  In such a case the colour of an edge of length $l$ may be written $c(l)$, or $c(l,U)$ where necessary to avoid ambiguity.  

A colouring of $U$ is {\it cyclic} if and only if (a) it is linear, and (b) $c(l)=c(m-l)$ for all $l$ such that $1 \le l \le m-1$.  

A colouring of $U$ is {\it repeating} in a colour $s$ with period $\pi$ if and only if (a) it is linear, and (b) $c(l)=s$ implies $c(l)=c(l + \pi)$ for all $l$ such that $1 \le l \le m-1-\pi$.  

The {\it clique number} of graph $U$ in colour $s$ is the largest integer $i$ such that $U$ contains a subgraph which is a copy of $K_i$ in that colour.

\section{Generalised Linear Graph Construction}

We shall set out general conditions for compounding two linear prototypes using a {\it triangle-free template}, or {\it tf-template}, in a relatively simple construction process, which includes that described in \cite{Rowley1} as a special case.

We first define a suitable form for the template, which denotes a subset of the first prototype with specified properties.  We go on to describe the construction and colouring of a graph $W$ in terms of the sets of its edge-lengths, as derived from the prototypes. For each of the colours of the first prototype, we prove that the clique number of $W$, however large, does not exceed that of $W'$, where $W'$ is derived from the graph $W$ by removing vertices.  Then, for all the colours of the second prototype, we prove that the clique-number of $W$ is the same as for the prototype. 

The construction is thus shown to give rise to a family of graphs of essentially unlimited size, the properties of which can be tested by searching only a limited subgraph, covering only a subset of their colours.  We note that, once tested in this way, the first prototype can be re-used with any second prototype without further testing.  This property radically limits the computing resources needed for establishing lower bounds using this method.  

\begin{definition}
 \label{Def:Def-1}
If the subset $\Theta$ comprising all the edge-lengths of a linear Ramsey graph $U(k_1, {\dots} \, ,k_{q-1}, 3; m)$ with colour $q$ contains the edge-length $m-1$, it is a {\it tf-template} for $U$.  
\end{definition}

The subset and the corresponding induced subgraph of $U$ can be identified with each other without confusion.  It should be clear that this definition makes $\Theta$ triangle-free. 
 
\bigskip    
\bigskip

\begin{theorem}
  \label{C-Thm3}
(Generalised Construction Theorem)\\
Given linear Ramsey graphs $U(k_1, {\dots} \, ,k_{q-1}, 3; m)$ and $V(k_{q+1}, {\dots} \, ,k_{q+r}; n)$, where $U$ contains a tf-template $\Theta$, it is possible to construct a linear Ramsey graph $W$ of order $(m-1)(n-1)+1+\phi$, for some $\phi$, where $0 \le \phi < m-1$.  The colour of the template is eliminated during construction: and if the maximum clique number in $U$ is $Q$, then the clique number of $W$ in any of the remaining colours of $U$ is no greater than the clique number of that colour in the subgraph of $W$ induced by the first $Q.(m-2)$ edge-lengths of $W$.  The clique number of $W$ in any of the colours of $V$ is the same as that of $V$.  
\end{theorem}

In intuitive terms, the theorem depends on a relatively simple construction process, building repeated copies of $U$, except that the colour of the template within each copy varies according to the colours of $V$, and so the colour $q$ is eliminated.  A partial copy of $U$ of order $\phi$, which cannot include any part of the template, is added as a final step.   

\begin{proof}
We begin by considering the set of lengths of all the edges of $U$, which we call $L$, consisting of the integers $\{\, l \mid 1 \le l \le m-1 \,\}$.  We know that a linear colouring gives rise to a natural partition of that set into subsets $L_s$ containing the lengths of edges of each colour $s$. That is, for $1 \leq s \leq q \,$:

\hspace{80pt} $L_s = \{\, l \in L \mid c(l) = s \}.$

As proved by Giraud in \cite{Giraud0}, if such a linear graph $U$ contains a copy of $K_{k_s}$ in colour $s$, then there exists a subset of the set $L_s$ of order $k_s-1$ such that each of the members of the subset and all of their non-zero pairwise absolute differences are contained in $L_s$.  That paper is presented in French, but it is simple to demonstrate that result and its converse, following his reasoning.

We first select a copy of $K_{k_s}$ and then reduce the indices of its vertex set by a uniform amount, so that the modified vertex set includes $u_0$.  The differences $(u_i - u_0)$ constitute the necessary subset of $L_s$.  Conversely, if a subset exists as described, in $L_s$, one can construct a set of all the vertices $u_i \in U$ having index-numbers $i$ in the subset.  Taking the union of that set of vertices with $u_0$ gives us the vertices of a copy of $K_{k_s}$ in $U$.  This result provides the basis for our proof.

The construction builds new subsets $L''_s$ of the edge-lengths of $W$.

For $1 \leq s \leq q-1 \,$, this proceeds as follows:

\hspace{80pt} $L'_s = \{\, l + (\mu - 1)(m-1) \mid l \in L_s, 1 \le \mu \le n-1 \}$.

Now let $\phi$ be the largest integer such that $\phi < \inf(l \in L_q)$ and for $1 \leq s \leq q-1$, define

\hspace{80pt} $L''_s = L'_s \cup \{\, l + (n-1)(m-1) \mid l \in L_s, 1 \le l \le \phi \}$.

Then for $1 \leq s \leq r$,

\hspace{80pt} $L''_{q+s} = \{\, l + (\mu - 1)(m-1) \mid l \in L_q$, $1 \le \mu \le n-1$, $c(\mu, V) = q + s \}$.

It is simple to verify that the new graph $W$ is of the required order, is well-defined and linear; and that it is repeating with period $m-1$ in all the colours of $U$ except for $q$, which has been eliminated in the construction -- so that there is no $L'_q$ or $L''_q$. 

We wish to know if there is a subset of the edge-lengths of $W$ that leads to a monochromatic copy of $K_{k_{s}}$ in colour $s$, where $1 \le s \le q-1$.  We aim to prove that if there is such a subset, then the edge-lengths contained in it must be less than a known bound.  

Assume now that there is such a subset, containing edge-lengths $\{h_1, h_2, \dots, h_{k_s - 1}\}$, and assume without loss of generality that these are strictly increasing. 

We can see from the definition that if the edge-length $h_1$ is greater than $m-1$ we can reduce all the edge-lengths by $m-1$ without changing the colour of any edge-length, or any difference between any pair of edge-lengths.  If any consecutive pair $h_{i_1}, h_{i_2}$ of these edge-lengths has a difference greater than $m-1$, we can reduce all the edge-lengths from $h_{i_2}$ upwards by $m-1$ without changing the colour of those edge-lengths, or any difference between any pair of edge-lengths.  As a result of repeating these steps, we can ensure that $h_{k_s - 1}$ does not exceed $(k_s - 1)(m-1)$: and since we also know that the colour of an edge-length which is any multiple of $(m-1)$ in $W$ is a colour derived from $V$, we can say that $h_{k_s - 1} \le (k_s - 1)(m-2)$.  

We have thus proved that when searching for potential copies of $K_{k_{s}}$ in the colour $s$ of $U$, it is not necessary to examine sets in $W$ containing edge-lengths in excess of $(k_s - 1)(m-2)$.  So if $Q = \sup\{k_s-1 \mid 1 \le s \le q-1 \}$, we have shown that the clique number of $W$ in any of the colours of $U$ is no greater than the clique number of that colour in the subgraph of $W$ induced by the first $Q.(m-2)$ edge-lengths of $W$, as required.  

We now consider the colours of $V$.  Assume that for any s, there are two edge-lengths $\mu_1$ and $\mu_2$ in $V$, both of colour $q+s$.  Assume that $\mu_1 < \mu_2$, and that their difference is also of that colour.  It follows directly from the definitions above that the edge-lengths $\mu_1(m-1)$, $\mu_2(m-1)$ and $(\mu_2 - \mu_1)(m-1)$ are all members of the subset $L''_{q+s}$.  Thus if there is a triangle of that colour in $V$, there is a triangle of the same colour in $L''_{q+s}$, and hence the same for any larger $K_n$.  

It is slightly more complex to prove the converse.  We assume there are two edge-lengths in $L''_{q+s}$, $l_i = h_i +(\mu_i-1)(m-1)$ for $i = 1, 2$, where $1 \le \mu_i \le n-1$ and $1 \le h_i \le m-1$.  Obviously they are both of colour $q+s$, and both $\mu_i$ are of the same colour in V.  We are interested in the possible colours of $\mid l_2 - l_1 \mid$ in $W$.  

Assume without loss of generality that $h_1 \le h_2$.  We note that, because the template $\Theta$ is triangle-free, the colour of $h_2 - h_1$ is not a colour of $V$.

If $\mu_2 \ge \mu_1$, then $l_2 \ge l_1$, and the colour of $l_2 - l_1$ is the same as the colour of $h_2 - h_1$, so it is not a colour of $V$.  Therefore $l_2 - l_1$ is not a member of $L''_{q+s}$.    

Assume now that $\mu_2 < \mu_1$.  The difference $l_1 - l_2$ is now positive and equal (after some rearrangement) to $((\mu_1 - \mu_2)-1)(m-1) + ((m-1) - (h_2 - h_1))$.  It is clear that $1 \le h_2 - h_1 \le m-2$, and if $c(l_1 - l_2) = q+s$, it follows that the colour of $(\mu_1 - \mu_2)$ in $V$ must also be $q+s$.  So if there is such a triangle in $W$, then there must be a triangle of the same colour in $V$, since it contains the edge-lengths $\mu_1$, $\mu_2$ and their difference. Again, it follows easily that the same applies for any larger $K_n$.

This completes the proof that the clique number of $W$ in any colour of $V$ is the same as the clique number of $V$ in that colour, which proves the theorem.

\end{proof}

This very general result has also been of some practical use in improving lower bounds for some significant `small' multicolour Ramsey numbers.  In addition it has improved our knowledge of the orders of many larger linear Ramsey graphs.  

It has been observed that there are many useful cases where the upper bound $Q.(m-2)$ is not tight, and also cases that can be found by simple searches where the clique number of $W$ in any of the colours of $U$ is the same as in $U$.  

It is also observed that when the subgraph of $U$ induced by the first $\phi$ edge-lengths is cyclic; when the colours of the last $(m-1)-\phi$ edge-lengths of $U$ are also symmetrically arranged (by reflection); and when $V$ is cyclic; then the constructed graph $W$ is also cyclic.  These conditions are illustrated pictorially by Figure 2 below.

Lastly, we note that if the prototype graph $U$ is taken to have even order, say $2p$, and the template is taken as the set $\{ p, p+1, p+2, \dots , 2p-1\}$, then the construction is identical to the construction in \cite{Rowley1} and $\phi = p-1$.

\section{New Lower Bounds for Several $R_r(k)$}
The method above has been used to construct a cyclic $(6,6,6; 1105)$-graph and a\\ $(5,5,5,5; 4072)$-graph.  For the first construction, the graph $U$ was taken as a very simple $(6,3; 12)$-graph and for the second a closely related $(5,3; 10)$-graph.  In the first case, the graph $V$ was taken as the $(6,6; 101)$-graph created by Kalbfleisch \cite{Kalb}, and in the second case an unpublished $(5,5,5; 453)$-graph derived by Exoo, which supersedes the results in \cite{XXER}.  Each graph $U$ has a tf-template as defined, and the constructed graphs have been verified completely, without relying on the theorem. These two values improve on the best lower bounds quoted in \cite{RadzDS}.  

The graphs $U$ can be represented in the following manner, which shows the colour to which each edge-length is mapped.  The tf-templates in Figure 1 are in colour 2 (blue).   

\begin{figure}[!ht]
  \begin{center}
    \includegraphics{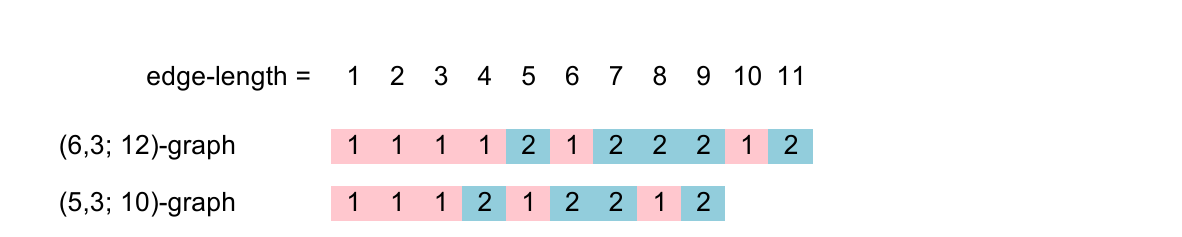}
  \end{center}
  \caption{\label{Fig:01} Graphs {\it U} used as first prototypes.}
\end{figure}

An intuitive picture of the compound colouring can be gained by illustrating the colours for each edge-length in a two-dimensional table.  Figure 2 shows the result of compounding the first template above with the unique $(3,4; 8)$-graph, purely to show the pattern of the colours.  The resulting $(3,4,6; 82)$-graph is useful as an illustration of the patterns involved, but not remarkable in itself.  

\begin{figure}[!ht]
  \begin{center}
    \includegraphics{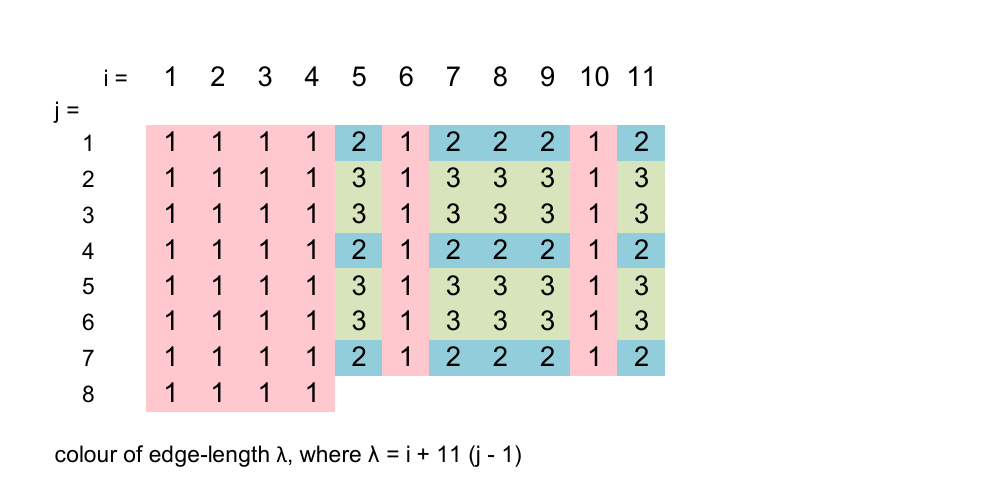}
  \end{center}
  \caption{\label{Fig:02} Illustration of a compound colouring.}
\end{figure}

\FloatBarrier

Table 1 shows the orders of a range of compound graphs constructed using first prototypes with two or three colours.  Note that the ` \& ' symbol denotes the compounding described in this paper, and ` * ' denotes the construction in \cite{Rowley1}.  

It is noted that certain values in Table 1 are highlighted to show improvements that exceed the values in the 2017 edition of the Radziszowski Dynamic Survey \cite{RadzDS}.  Many of these values are now included in the 2021 edition. 

\begin{table}[!ht]
  \begin{center}
    \includegraphics{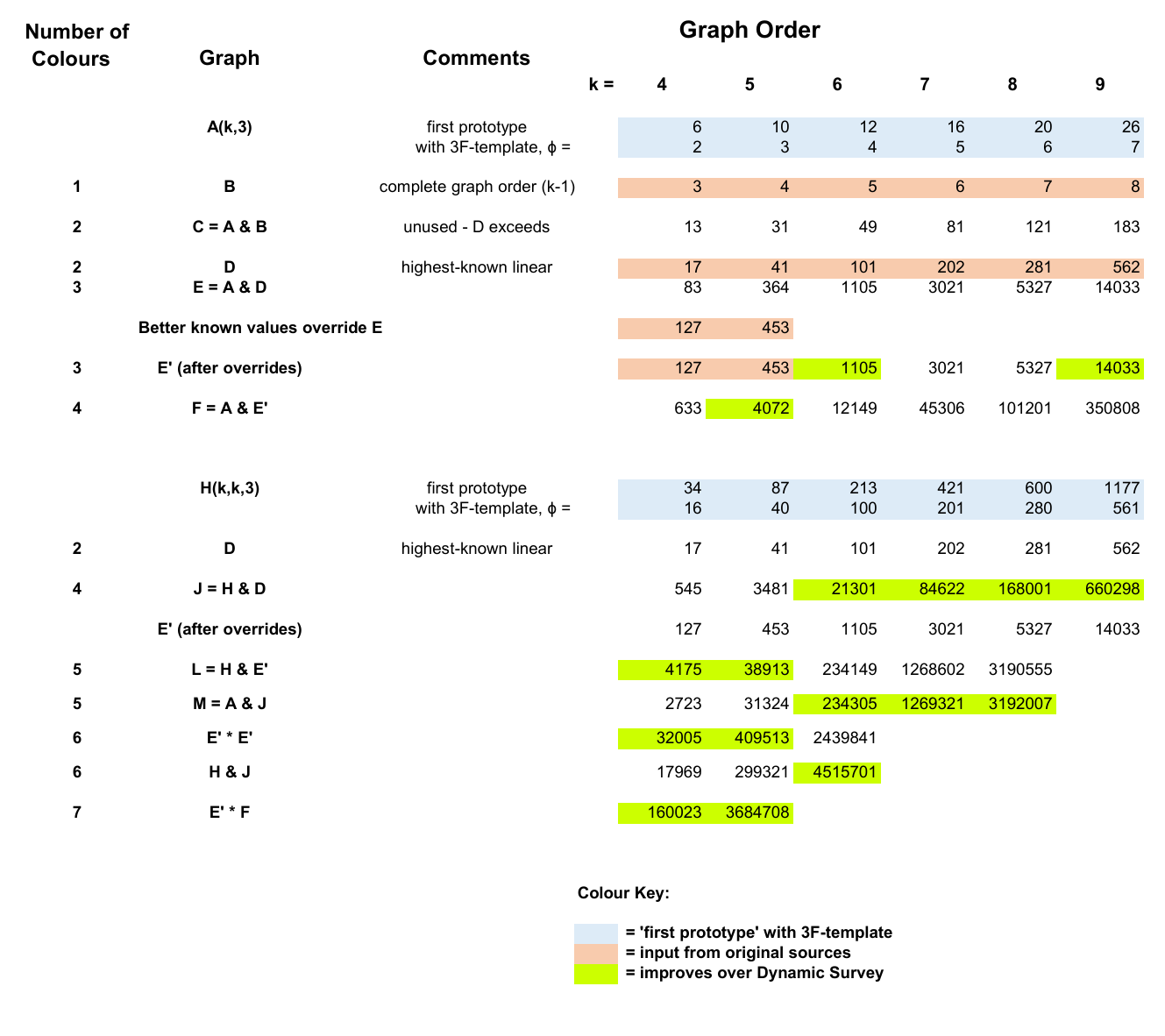}
  \end{center}
  \caption{\label{Table:01} Calculations of graph orders to support lower bounds for $R_r(k)$.}
\end{table}

\FloatBarrier

\section{Updated Tables of Results}

The tables below update the results of the previous paper \cite{Rowley1}, allowing for the inclusion of many new graphs constructed using the method mentioned above.  They show the revised linear graph orders and associated `growth factors'.  These growth factors represent the limiting growth rates of an infinite series of graphs which can be produced by the construction set out in \cite{Rowley1}, by repeatedly using the graph concerned in the compounding process.  Numbers revised since the publication of \cite{Rowley2} are shown in green or blue. Orders above 10 million have been excluded.  Bold text indicates numbers exceeding those shown in the 2017 Radziszowski Dynamic Survey \cite{RadzDS}, although some have been included in the 2021 edition. 

All these graphs have been computer-tested to the extent implied by the theorem, and in some cases fully, up to 3000 vertices or more.  

\begin{table}[!ht]
  \begin{center}
    \includegraphics{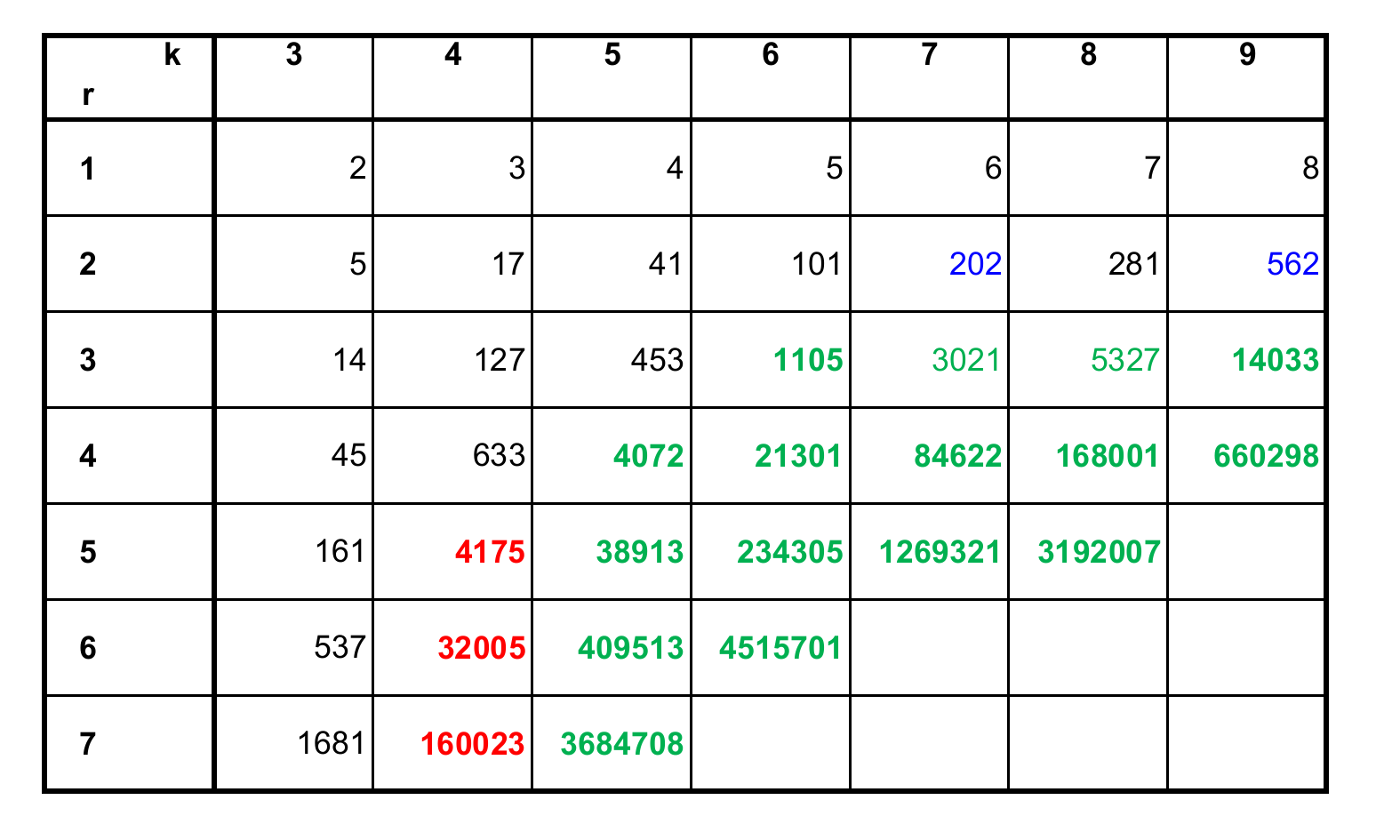}
  \end{center}
  \caption{\label{Table:02} Highest order of linear Ramsey graphs $U_r(k)$ known to the author.}
\end{table}

\medskip

\FloatBarrier

The factors $g_r(k)$ in Table 3 are defined as $(2m-1)^{1/r}$ where m is the order of the relevant graph in Table 2.  Those highlighted in bold indicate are most effective as lower bounds for $\Gamma(k)$.  

\begin{table}[!ht]
  \begin{center}
    \includegraphics{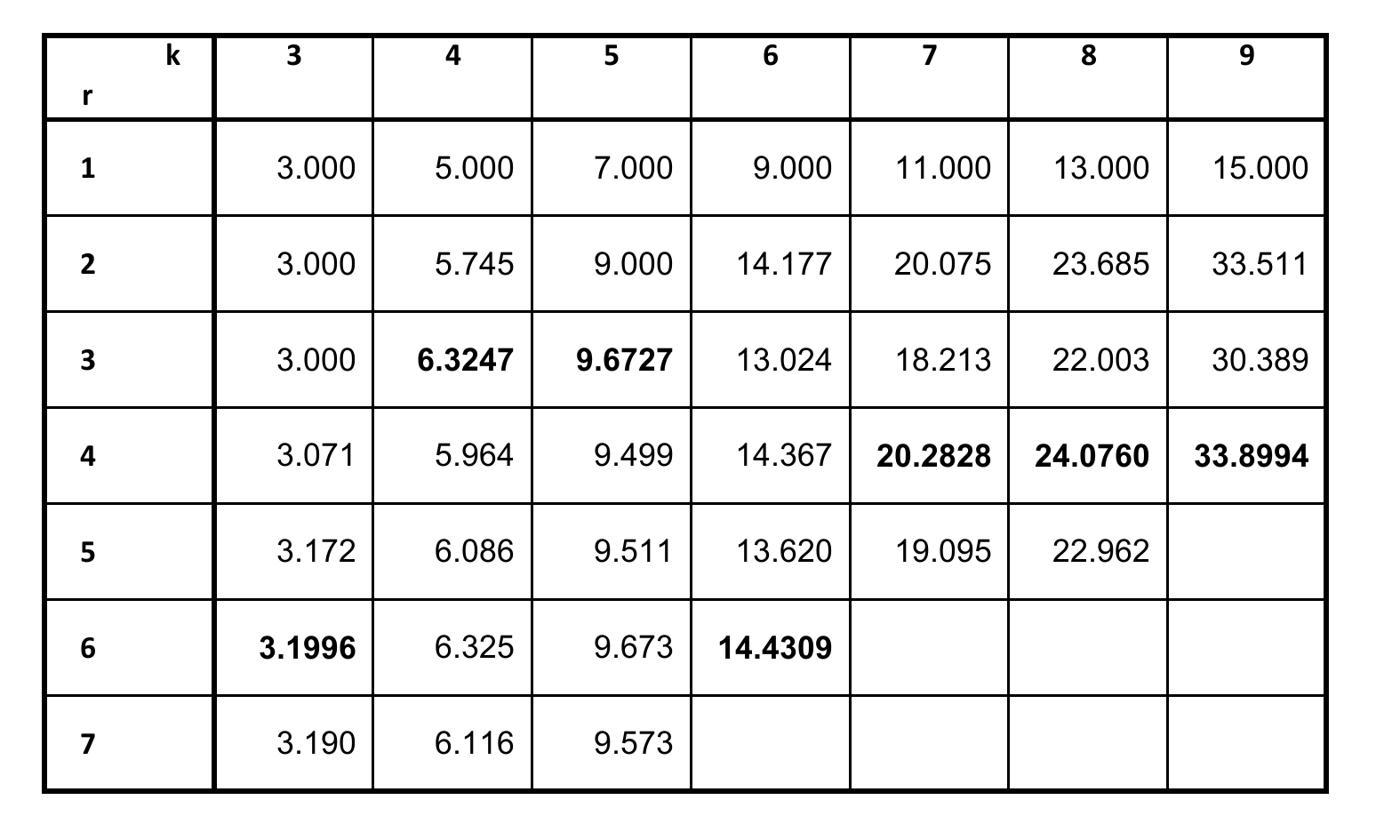}
  \end{center}
  \caption{\label{Table:03} Growth factors $g_r(k)$ calculated from the data in Table 2.}
\end{table}

\FloatBarrier

In fact, the lower bounds on $\Gamma(k)$ derived from Table 3 can be improved somewhat by the application of the construction in this paper and some simple arithmetic.  For example, it can be shown that the sequence of increasing values for $k = 6$, and even $r$, leads into an infinite sequence of graphs whose growth rate converges to the square root of $m-1$, where $m$ is the order of the original first prototype (graph $H$ in Table 1).  In this way, we see that $\Gamma(6) \ge \sqrt{212} = 14.560\dots$, $\Gamma(7) \ge \sqrt{420} = 20.493\dots$, $\Gamma(8) \ge \sqrt{599} = 24.474\dots$ and $\Gamma(9) \ge \sqrt{1176} = 34.292\dots$, etc.  

The same reasoning, given the existence of a favourable $(3,3,3,3,3,3; 377)$-graph, shows that $\Gamma(3) \ge \sqrt[5]{376} = 3.273\dots$\,.  Although this graph produces no corresponding improvements in lower bounds for Ramsey numbers $R_r(3)$ below $R_8(3)$, it is noted that the techniques described in this paper enable construction of graphs establishing that $R_8(3) \ge 5288$, $R_9(3) \ge 17696$, $R_{10}(3) \ge 60322$, $R_{11}(3) \ge  201698$ and $R_{12}(3) \ge 631842$.

Rearranging, we see that the same graph also provides an improved lower bound for $\Gamma(3)$, since it demonstrates that $R{_r}(3) > c.(3.273)^r$, for some fixed $c$ and for $r$ sufficiently large.

An interesting special case depends on a formula derived in \cite{Abbott} and mentioned in \cite{AbbH} -- namely that $R_r(5) \ge (R_r(3) - 1)^2 + 1$.  It follows that $\Gamma(5) \ge 10.717\dots$.   

\bigskip    

\section{Some Conclusions}

Constructive approaches to graph colouring based on compounding linear graphs have the favourable property that they restrict the need to search for occurrences of $K_n$ within the compound graph.  There is a trade-off between this favourable property and the complexity of the prototypes that can be used in the compounding process.  

The generalised construction above expands the range of prototypes that can be used in the compounding, at the cost of increasing the search requirements.  Even so, the resulting graphs mentioned here can be tested in quite manageable times.  Although there is an added need to search for suitable prototypes containing tf-templates, it is useful in practice that prototypes with tf-templates can be freely used in conjunction with any other linear prototype graph. 

A great strength of the current construction -- like that in \cite{Rowley1} -- is that it can be applied to any linear prototype graphs with the desired Ramsey properties, irrespective of their clique numbers (with the notable exception of the tf-template).

Its power is demonstrated by the fact that many values of $R_r(k)$, and the corresponding lower bounds on $\Gamma(k)$, can be improved using this construction.  This paper contains a range of examples. This improvement arises because the previous construction cannot increase the limiting growth rate of any series of derived graphs beyond the maximum growth factor for its prototypes, since it produces a weighted mean of their factors.  In contrast, the current construction gains an advantage by eliminating the template colour, so that in many cases, further improved bounds on $\Gamma(k)$ can be derived as mentioned in section 5, by taking roots.  

In a particular case noted above, the current construction provides an improved lower bound for the limiting growth rate of $R{_r}(3)$, by demonstrating that $R{_r}(3) > c.(3.273)^r$, for some fixed $c$ and sufficiently large $r$, which is equally relevant to the corresponding Schur numbers.

The most obvious exception to these improvements is the case of $k = 4$, where the current construction has not yet yielded any improvement in the smaller graph orders.  It may be that this results from the high criticality of the $(4,4; 17)$- and $(4,4,4; 127)$-graphs, which seems to inhibit the construction of large first prototypes with useful templates.

Finally, it is interesting to consider the improvements to the lower bounds on $R_3(6)$ to $1106$, and on $R_3(9)$ to $14034$, relative to non-linear graphs.  The previous value of $1070$ on $R_3(6)$ was established by Mathon in \cite{Mathon} based on a cyclic graph. By contrast, the lower bounds on $R_3(k)$ in the Dynamic Survey were established using non-linear graphs for $k > 6$.  Although the order of the new linear $R_3(9)$-graph exceeds the previous best lower bound, the existing bounds for $k = 7, 8$ remain above those achieved so far using linear constructions.  Indeed, the existence of the graph $H(8,8,3)$ in Table 1 shows the existence of a graph with order 880, which (by `quadrupling' twice in succession to produce a non-linear graph) shows that $R_3(9) \ge 14081$. 

\subsection*{Acknowledgements}
I once again record my warmest thanks to my wife Joan for her continued support for this work, and to my son William for sharing his ideas for improved search programs. 

Many thanks also to Geoffrey Exoo, for making available the remarkable\\ (5,5,5; 453)-graph mentioned in the paper, and to Stanis{\l}aw Radziszowski, for kindly bringing it my attention.  

Finally, my thanks also go to the anonymous reviewers, whose comments considerably strengthened this paper.


\begin{thebibliography}{12}

\bibitem{Abbott} H.L.~Abbott, \newblock Some Problems in Combinatorial Analysis, \newblock {\em Ph.D. Thesis},\\  Department of Mathematical and Statistical Sciences, University of Alberta, \\ Edmonton (1965).

\bibitem{AbbH} H.L.~Abbott and D.~Hanson, \newblock A Problem of Schur and its Generalizations, \\ \newblock {\em Acta Arithmetica}, {\bf \small XX} (1972), 175--187.

\bibitem{AbbM} H.L.~Abbott and L.~Moser, \newblock Sum-free Sets of Integers, \\ \newblock {\em Acta Arithmetica}, {\bf \small XI} (1966), 393--396.

\bibitem{Giraud0} Guy~R.~Giraud, \newblock Une g\'en\'eralisation des nombres et de l'in\'egalit\'e de Schur, \newblock {\em C.R. Acad. Sc. Paris}, S\'erie A, {\bf 266} (1968), 437--440.

\bibitem{Giraud1} Guy~R.~Giraud, \newblock Minorations de certains nombres de Ramsey binaires \\ par les nombres de Schur g\'en\'eralis\'es, \newblock {\em C.R. Acad. Sc. Paris}, S\'erie A, {\bf 266} (1968), 481--483.

\bibitem{Kalb} J.G.~Kalbfleisch, \newblock Chromatic Graphs and Ramsey's Theorem, \newblock {\em Ph.D. Thesis}, University of Waterloo, Canada (1966).

\bibitem{Mathon} R.~Mathon, \newblock Lower Bounds for Ramsey Numbers and Association Schemes, \\ \newblock {\em Journal of Combinatorial Theory}, Series B, {\bf 42}, (1987), 122--127.

\bibitem{RadzDS} S.P.~Radziszowski, \newblock Small Ramsey Numbers, \newblock {\em Electron. J. Combin.} (2017), {\#}DS1.

\bibitem{Rowley1} F.~Rowley, \newblock Constructive Lower Bounds for Ramsey Numbers from Linear Graphs, \newblock {\em Australasian J. Combin.} {\bf 68(3)} (2017), 385--395.  

\bibitem{Rowley2} F.~Rowley, \newblock Some Further Results in Ramsey Graph Construction, \newblock {\em Australasian J. Combin.} {\bf 78(1)} (2020), 1--10. 

\bibitem{Sch1} I.~Schur, \"{U}ber die Kongruenz $x^m + y^m \equiv z^m~(mod~p)$, {\it Jahresber. Dtsch. Math.-Ver.}  {\bf 25} (1917), 114-116. ~~(\url{http://eudml.org/doc/145475}).

\bibitem{XXER} X.~Xu, Z.~Xie, G.~Exoo and S.P.~Radziszowski, \newblock Constructive Lower Bounds on Classical Multicolour Ramsey Numbers, \newblock {\em Electron. J. Combin.} {\bf 11} (2004), {\#}R35.

\end{thebibliography}

\end{document}